\documentclass{article}
\usepackage{latexsym,bm}
\usepackage{amsmath,amsfonts,amssymb}

\DeclareMathOperator{\ch}{char} \DeclareMathOperator{\soc}{soc}
\DeclareMathOperator{\Hom}{Hom} \DeclareMathOperator{\res}{res}
\DeclareMathOperator{\bH}{H} \DeclareMathOperator{\id}{id}
 
\DeclareMathOperator{\rad}{rad}

\newcommand{\ksl}{\widehat{\mathfrak{sl}}}
\renewcommand{\l}{\ell}
\renewcommand{\H}{\mathcal{H}}
\newcommand{\lam}{\lambda}
\newcommand{\BS}{\mathfrak S}
\newcommand{\ts}{\widetilde S}
\newcommand{\td}{\widetilde D}
\newcommand{\Z}{\mathbb Z}
\newcommand{\qed}{\hfill\text{$\square$}}

\newtheorem{thm}{Theorem}\newtheorem{cor}{Corollary}
\newtheorem{lem}{Lemma}\newtheorem{defs}{Definition}
\newtheorem{example}{Example}
\numberwithin{example}{section} \numberwithin{equation}{section}
\numberwithin{prop}{section} \numberwithin{thm}{section}
\numberwithin{lem}{section} \numberwithin{defs}{section}
\numberwithin{cor}{section}

\title{Branching rules for Hecke Algebras of Type $D_{n}$\thanks{Research support by Alexander von
Humboldt Foundation, National Natural Science Foundation of
China (Project 10401005) and 21st Century COE program. Part of the
work was done while the author was visiting the Research Institute
for Mathematical Science at Kyoto University. The author thanks
the hospitality of RIMS during the writing of this paper.}}
\author{Jun Hu$^{\star}$ \\[5pt]
$^{\star}$Department of Applied Mathematics\\
Beijing Institute of Technology\\
Beijing, 100081, P.R. China\\[1.5pt]
E-mail: junhu303@yahoo.com.cn}

\date{}
\begin{document}
\maketitle

\begin{abstract} In
this paper we study the branching problems for Hecke algebra
$\H(D_n)$ of type $D_n$. We explicitly describe the decompositions
of the socle of the restriction of each irreducible
$\H(D_n)$-representation to $\H(D_{n-1})$ into irreducible modules
by using the corresponding results for type $B$ Hecke algebras. In
particular, we show that any such restrictions are always
multiplicity free.
\end{abstract}

\section{Preliminaries}

Let $W(B_{n})$ be the Weyl group of type $B_{n}$. It is a finite
group with generators $\{s_{0}, s_{1},\cdots,s_{n-1}\}$ and
relations $$ \aligned
&s_i^2=1,\quad\text{for\, $0\leq i\leq n-1$,}\\
&s_0s_1s_0s_1=s_1s_0s_1s_0,\\
&s_is_{i+1}s_i=s_{i+1}s_{i}s_{i+1},\quad\text{for\, $1\leq i\leq n-2$,}\\
&s_is_j=s_js_i,\quad\text{for\, $0\leq i<j-1\leq n-2$.}\endaligned
$$
Let $u:=s_0s_1s_0$. The subgroup of $W(B_{n})$ generated by $\{u,
s_{1},\cdots,s_{n-1}\}$ is a Weyl group of type $D_{n}$. We denote
it by $W(D_{n})$. It has a presentation with generators $\{u,
s_{1},\cdots,s_{n-1}\}$ and relations $$ \aligned
&u^2=1=s_i^2,\quad\text{for\, $1\leq i\leq n-1$,}\\
&us_2u=s_2us_2,\quad s_is_{i+1}s_i=s_{i+1}s_{i}s_{i+1},\quad\text{for\, $1\leq i\leq n-2$,}\\
&us_1=s_1u,\quad us_i=s_iu,\quad\text{for \,$3\leq i\leq n-1$,}\\
&s_is_j=s_js_i,\quad\text{for\, $1\leq i<j-1\leq n-2$.}\endaligned
$$
Let $\BS_n$ be the symmetric group on $n$ letters. It is
well-known that $W(B_{n})\cong (\Z/2\Z)^{n}\rtimes{\BS}_{n},\,
W(D_{n})\cong (\Z/2\Z)^{n-1}\rtimes{\BS}_{n}$, and the subgroup
generated by $s_{1}, s_{2}, \cdots$ $,s_{n-1}$ (respectively $u,
s_{2}, \cdots, s_{n-1}$) can be identified with the symmetric
group $\BS_n$.

Let $K$ be a field. {\it Throughout this paper we assume that
$\ch{K}\neq 2$}. Let $q, Q$ be two invertible elements in $K$.
There is a Hecke algebra $\H_{q, Q}(B_{n})$ with parameters $q, Q$
associated to $W(B_{n})$ (see \cite{DJ}). In this paper we will
only be concerned with the special case where $Q=1$, i.e.,
$\H(B_n):=\H_{q, 1}(B_{n})$. By definition, $\H(B_n)$ is an
associative algebra with generators $T_0, T_1,\cdots, T_{n-1}$ and
relations
$$\aligned
&T_0^2=1,\quad (T_i+1)(T_i-q)=0,\quad\text{for\, $1\leq i\leq n-1$,}\\
&T_0T_1T_0T_1=T_1T_0T_1T_0,\\
&T_iT_{i+1}T_i=T_{i+1}T_{i}T_{i+1},\quad\text{for\, $1\leq i\leq n-2$,}\\
&T_iT_j=T_jT_i,\quad\text{for\, $0\leq i<j-1\leq n-2$.}\endaligned
$$
Let $T_u:=T_0T_1T_0$. The subalgebra of $\H(B_{n})$ generated by
$\{T_u,T_{1},\cdots,T_{n-1}\}$ is isomorphic to a Hecke algebra of
type $D_n$, i.e., the Hecke algebra associated to the Weyl group
$W(D_{n})$. We denote it by $\H(D_{n})$. It has a presentation
with generators $\{T_u, T_{1},\cdots, T_{n-1}\}$ and relations
$$\aligned
&(T_u+1)(T_u-q)=0,\quad (T_i+1)(T_i-q)=0,\quad\text{for\, $1\leq i\leq n-1$,}\\
&T_uT_{2}T_u=T_{2}T_{u}T_{2},\quad
T_iT_{i+1}T_i=T_{i+1}T_{i}T_{i+1},
\quad\text{for\, $1\leq i\leq n-2$,}\\
&T_uT_1=T_1T_u,\quad T_uT_i=T_iT_u,\quad\text{for\, $3\leq i\leq n-1$,}\\
&T_iT_j=T_jT_i,\quad\text{for\, $1\leq i<j-1\leq n-2$,}
\endaligned
$$

By \cite{DJMu}, $\H(B_n)$ is a cellular algebra in the sense of
\cite{GL}. For each bipartition $\lam=(\lam^{(1)},\lam^{(2)})$ of
$n$, there is a {\it Specht module} $\ts^{\lam}$, and a naturally
defined bilinear form on $\ts^{\lam}$. Let $\td^{\lam}$ be the
quotient of $\ts^{\lam}$ modulo the radical of that form. We have
that

\begin{lem} [\rm \cite{DJMu}] 1) Every simple $\H(B_{n})$ module is a composition factor of some
$\ts^{\lam}$. When $\H(B_{n})$ is semi-simple, each $\ts^{\lam}$
is absolutely irreducible and they form a complete set of pairwise
non-isomorphic simple $\H(B_{n})$-modules.

2) If $\td^{\mu}\neq 0$ is a composition factor of $\ts^{\lambda}$
then $\lam\trianglerighteq\mu$, and every composition factor of
$\ts^{\lam}$ is isomorphic to some $\td^{\mu}$ with
$\lam\trianglerighteq\mu$, where $\trianglerighteq$ is the dominance order defined in \cite{DJMu}. 
If $\td^{\lam}\neq 0$ then the
composition multiplicity of $\td^{\lam}$ in $\ts^{\lam}$ is one.

3) The set $\Bigl\{\td^{\lam}\Bigm|\text{$\lam$ is a bipartition
of $n$ and $\td^{\lam}\neq 0$}\Bigr\}$ forms a complete set of
pairwise non-isomorphic simple $\H(B_{n})$-modules.

4) $\H(B_{n})$ is semi-simple if and only if $\,\,
2\biggl(\prod_{i=1}^{n-1}\bigl(1+q^i\bigr)\biggr)\biggl(\prod_{i=1}^{n}\bigl(1+q+q^2+
\cdots+q^{i-1}\bigr)\biggr)\neq 0$. In that case, it is also split
semi-simple.
\end{lem}

Let $\tau$ be the involutive $K$-algebra automorphism of
$\H(B_{n})$ which is defined on generators by
$\tau(T_{1})=T_{0}T_{1}T_{0}, \tau(T_i)=T_i,\,\,\forall\, i\neq
1$. Then $\tau$ maps $\H(D_{n})$ isomorphically onto $\H(D_{n})$.
Let $\sigma$ be the involutive $K$-algebra automorphism of
$\H(B_{n})$ which is defined on generators by
$\sigma(T_{0})=-T_{0}, \sigma(T_i)=T_i,\,\,\forall\, i\neq 0$.
Then $\sigma\downarrow_{\H(D_n)}=\id$. Let $\widetilde{\mathcal
P}_n$ be the set of all the bipartitions of $n$.

\addtocounter{defs}{1}
\begin{defs} For each $\lam=(\lam^{(1)},\lam^{(2)})\in\widetilde{\mathcal{P}}_n$, we
define
$\widehat{\lam}=(\lam^{(2)},\lam^{(1)})\in\widetilde{\mathcal{P}}_n$.
\end{defs}

Let ${\mathcal
P}_n:=\bigl\{\lam\in\widetilde{\mathcal{P}}_n\bigm|\td^{\lam}\neq
0\bigr\}$. There is an involution $\bH$ defined on ${\mathcal
P}_n$ such that
$\Bigl(\td^{\lam}\Bigr)^{\sigma}\cong\td^{\bH(\lam)}$ for any
$\lam\in{\mathcal P}_n$. In particular,
$\td^{\lam}\!\!\downarrow_{\H(D_n)}\cong\td^{\bH(\lam)}\!\!\downarrow_{\H(D_n)}$.
We define an equivalence relation on ${\mathcal P}_n$ by
$\lam\approx\mu$ if and only if $\mu=\bH(\lam)$. By the results in
\cite{P}, \cite{Hu1}, \cite{Hu2} and \cite{Hu3}, we have that

\addtocounter{lem}{1}
\begin{lem} [\rm \cite{Hu1}] Suppose that $\ch K\neq
2$ and $\H(D_n)$ is split over $K$. If $\lam\neq\bH(\lam)$, then
$\td^{\lam}\!\!\downarrow_{\H(D_n)}$ is irreducible; if
$\lam=\bH(\lam)$, then $\td^{\lam}\!\!\downarrow_{\H(D_n)}$ splits
into a direct sum of two $\H(D_n)$-submodules, say $D_{+}^{\lam}$
and $D_{-}^{\lam}$. Moreover, the set $$
\Bigl\{\td^{\lam}\!\!\downarrow_{\H(D_n)}\Bigm|\lam\in{\mathcal
P}_n/{\approx},\,\,\lam\neq\bH(\lam)\Bigr\}\bigcup
\Bigl\{D_{+}^{\lam}, D_{-}^{\lam}\Bigm|\lam\in{\mathcal
P}_n/{\approx},\,\,\lam=\bH(\lam)\Bigr\} $$ forms a complete set
of pairwise non-isomorphic irreducible $\H(D_n)$-modules.
\end{lem}

Let $e>1$ be a fixed integer. Let $\lam$ be a bipartition 
of $n$. For each node $\gamma=(i,j)$ of $\lam$, we define the residue of $\gamma$
to be $j-i+e\mathbb{Z}\in\mathbb{Z}/e\mathbb{Z}$. Then we have  
the notion of $e$-good (removable) nodes of $\lam$
(see \cite{AM} and \cite{Hu2}). For each integer $m\in\mathbb{N}$,
the set $\mathcal{K}_m$ of Kleshchev bipartitions of $m$ with
respect to $(\sqrt[e]{1};1,-1)$ is defined inductively by
\smallskip

(1)
$\mathcal{K}_0:=\Bigl\{\underline{\emptyset}:=\bigl(\emptyset,\emptyset\bigl)\Bigr\}$;

(2)
$\mathcal{K}_{m}:=\Bigl\{\lam\in\widetilde{\mathcal{P}}_{m}\Bigm|\begin{matrix}\text{$\lam$
is obtained from some $\mu\in\mathcal{K}_{m-1}$ by}\\
\text{adding an $e$-good node}\end{matrix}\Bigr\}$.
\medskip

The {\it Kleshchev's good lattice} with respect to
$(\sqrt[e]{1};1,-1)$ is the infinite graph whose vertices are the
Kleshchev bipartitions with respect to
$(\sqrt[e]{1};1,-1)$ and whose arrows are given by $$
\text{$\mu\overset{x}{\rightarrow}\lam$\quad$\Longleftrightarrow$\quad
$\lam$ is obtained from $\mu$ by adding an $e$-good $x$-node}.
$$ By a result of S. Ariki (see \cite{A}),
$\mathcal{P}_n=\mathcal{K}_n$ when $e$ is the smallest positive integer satisfying
$1+q+q^2+\cdots+q^{e-1}=0$.

\begin{lem} [\rm \cite{P}, \cite{Hu1}, \cite{Hu2}, \cite{Hu3}]

1) If $2\prod_{i=1}^{n-1} \bigl(1+q^i\bigr)\neq 0$ in $K$, then
$\H(D_n)$ is split over $K$. In this case, for each
$\lam=(\lam^{(1)},\lam^{(2)})\in\mathcal{P}_n$, we have that
$\bH(\lam)=\widehat{\lam}$. In particular,
$\td^{\lam}\!\!\downarrow_{\H(D_{n})}\cong\td^{\widehat\lam}\!\!\downarrow_{\H(D_{n})}$,\smallskip

2) If $\,\,\prod_{i=1}^{n-1} \bigl(1+q^i\bigr)=0$, $\,\ch K\neq 2$
and $\H(D_n)$ is split over $K$, then $q$ is a primitive $2\l$-th
root of unity in $K$ for some integer $1\leq\l<n$. In this case,
$\bH$ can be described as follows: if $\lam\in\mathcal{P}_n$ is a
Kleshchev bipartition with respect to $(\sqrt[2\l]{1};1,-1)$, and
$$
\underline{\emptyset}\stackrel{i_1}{\longrightarrow}\cdot\overset{i_2}{\longrightarrow}\cdot
\cdots\stackrel{i_n}{\longrightarrow}\lam, $$ is a path from
$\underline{\emptyset}:=(\emptyset,\emptyset)$ to $\lam$ in
Kleshchev's good lattice with respect to $(\sqrt[2\l]{1};1,-1)$.
Then, the sequence $$
\underline{\emptyset}\stackrel{i_1+\l}{\longrightarrow}\cdot\overset{i_2+\l}{\longrightarrow}\cdot
\cdots\stackrel{i_n+\l}{\longrightarrow}\bH(\lam), $$  also
defines a path in Kleshchev's good lattice with respect to
$(\sqrt[2\l]{1};1,-1)$, and it connects $\underline{\emptyset}$ to
$\bH(\lam)$.\footnote{This result was proved in \cite{Hu2} only in
the case where $K=\mathbb{C}$. For more general $K$, see appendix
of this paper.}

3) If $\,\ch K\neq 2$, then  $\H(D_{n})$ is semi-simple if and
only if $$
\biggl(\prod_{i=1}^{n-1}\bigl(1+q^i\bigr)\biggr)\biggl(\prod_{i=1}^{n}\bigl(1+q+q^2+
\cdots+q^{i-1}\bigr)\biggr)\neq 0.$$ In that case, it is also
split semi-simple.
\end{lem}

In this paper, we shall give the modular branching rule for
$\H(D_n)$. That is, for each irreducible $\H(D_n)$-module $D$, we
describe $\soc\bigl(D\!\!\!\downarrow_{\H(D_{n-1})}\bigr)$. The
discussion will be divided into two cases: the case where
$2\prod_{i=1}^{n-1}\bigl(1+q^i\bigr)\neq 0$ and the case where
$\prod_{i=1}^{n-1}\bigl(1+q^i\bigr)=0$ and $2\cdot 1_K\neq 0$. The
main results are presented in Theorem 2.5, Theorem 2.6, Corollary
2.8, Theorem 3.7, Theorem 3.9 and Corollary 3.11. It turns out that
our situation here bears much resemblance to the situation of
representations of the alternating group $A_n$ (which is a normal
subgroup in $\BS_n$ of index $2$), see \cite{B}, \cite{FK},
\cite{LLT} and \cite{Hu2}, and our results are largely motivated
by those in \cite{BO}, where the branching rules for the
representations of the alternating groups were deduced.
\medskip

Finally, in the appendix of this paper, we include a proof (which
is essentially due to Professor S. Ariki) of the fact that the
involution $\bH$ is independent of the base field $K$ as long as
$\ch K\neq 2$ and $\H(D_n)$ is split over $K$.

\smallskip\bigskip
\section{The case where
$2\prod_{i=1}^{n-1}\bigl(1+q^i\bigr)\neq 0$}

Let $\lam$ be a bipartition, $[\lam]$ be its Young diagram. To
simplify notation, we shall identify $\lam$ with $[\lam]$. Recall
that a {\it removable} node is a node of the boundary of $[\lam]$
which can be removed, while an {\it addable} node is a concave
corner on the rim of $[\lam]$ where a node can be added. {\it
Throughout this section, we shall assume that
$2\prod_{i=1}^{n-1}\bigl(1+q^i\bigr)\neq 0$ in $K$.} In
particular, Lemma 1.4(1) applies to both the Hecke algebra
$\H(D_n)$ and the Hecke algebra $\H(D_{n-1})$.

Let $l$ be the smallest positive integer $a$ such that
$1+q+q^2+\cdots+q^{a-1}=0$. If such an integer does not exist,
then we set $l=\infty$. A partition $\lam=(\lam_1,\lam_2,\cdots)$
is said to be $l$-restricted if $\lam_i-\lam_{i+1}<l$ for all $i$.

\begin{lem} [\rm \cite{DJ}, \cite{DR}] Suppose that
$\,2\prod_{i=1}^{n-1}\bigl(1+q^i\bigr)\neq 0$ in $K$. Then

1) for each bipartition $\lam=(\lam^{(1)},\lam^{(2)})$ of $n$, we
have that $\lam\in\mathcal{P}_n$ if and only if both $\lam^{(1)}$
and $\lam^{(2)}$ are $l$-restricted.

2) for each bipartition $\lam\in\mathcal{P}_n$, $\,\,
\soc\Bigl(\td^{\lam}\!\!\!\downarrow_{\H(B_{n-1})}\Bigr)\cong\bigoplus_{\mu\rightarrow\lam}
\td^{\mu}$, where $\mu\rightarrow\lam$ means that $\mu$ is
a bipartition of $\,\,n-1$ such that the Young diagram $[\mu]$ is obtained from the Young diagram 
$[\lam]$ by removing an $l$-good node.
\end{lem}

Note that in this case the Kleshchev's good lattice with respect to
$(\sqrt[l]{1};1,-1)$ is well-understandood. Namely, for any Kleshchev bipartition $\lam$, a removable node 
$\gamma$ is an $l$-good node of the bipartition $\lam$ if and only if $\gamma$ is an $l$-good node
of the partition $\lam^{(1)}$ or of the partition $\lam^{(2)}$. Here the notion of $l$-good nodes of 
partitions is defined in a similar way as $l$-good nodes of bipartitions, see \cite{Kl}, \cite{LLT} and 
\cite{AM} for details.\smallskip

We want to describe the decomposition of the socle of
$D\!\!\downarrow_{\H(D_{n-1})}$ into irreducible modules for each
irreducible $\H(D_n)$-module $D$. For each bipartition $\lam$ and
each removable node $A$ of $[\lam]$, we shall denote by
$\lam\setminus{A}$ the Young diagram (or equivalently,
bipartition) obtained by removing the node $A$ from $[\lam]$.

%The following definition is motivated by \cite{BO}.

\addtocounter{defs}{1}
\begin{defs} Let $\lam\in\mathcal{P}_n$.
Suppose that $\lam$ has an $l$-good node $A$ such that $$
\lam\setminus{A}=\widehat{(\lam\setminus{A})}=\widehat\lam\setminus{A'}
$$
for some removable node $A'$ in $\widehat\lam$. Then in this case
the node $A'$ is uniquely determined by $A$ and is also an $l$-good
node of $\widehat\lam$. We say that $\lam$ is almost symmetric and
$A'$ is the conjugate node of $A$.
\end{defs}

\addtocounter{example}{2}
\begin{example} Suppose that $n=5$, $l=\infty$.
$\lam=((2,1), (1^2)),\, \mu=((2), (1^3))$ are two bipartitions of
$5$. Let $A$ be the node which is in the first row and the second
column of the first component of $\lam$. Then $\lam$ is almost
symmetric with $\lam\setminus{A}=\widehat{(\lam\setminus{A})}$,
but $\mu$ is not almost symmetric.
\end{example}

\addtocounter{lem}{2}
\begin{lem} Let $\lam=(\lam^{(1)},\lam^{(2)})$ be a
bipartition of $n$. Suppose that $\lam$ is almost symmetric with
$\lam\setminus{A}=\widehat{(\lam\setminus{A})}$ for some removable
node $A$ of $\lam$. Then for any pairs of removable nodes $B, C$
of $\lam$ satisfying $C\neq A$, we have that $
\lam\setminus{B}\neq\widehat{(\lam\setminus{C})}$. In particular,
for any removable nodes $C$ of $\lam$ satisfying $C\neq A$, we
have that $$ \lam\setminus{C}\neq\widehat{(\lam\setminus{C})}. $$
\end{lem}

\noindent {Proof:} \, This is obvious. \hfill\qed

\addtocounter{thm}{4}
\begin{thm} Let $\lam\in\mathcal{P}_n$. Suppose that
$\lam$ is almost symmetric with
$\lam\setminus{A}=\widehat{(\lam\setminus{A})}$ for some $l$-good node
$A$ of $\lam$. Then $\lam\neq\widehat\lam$ and $$
\soc\Bigl(\td^{\lam}\!\!\downarrow_{\H(D_{n-1})}\Bigr)\cong
D_{+}^{\lam\setminus{A}}\bigoplus D_{-}^{\lam\setminus{A}}\bigoplus
\bigoplus_{\substack{C\in [\lam],\,\, C\neq A\\ \text{$C$ is
$l$-good}}}\td^{\lam\setminus{C}}\!\!\downarrow_{\H(D_{n-1})}.
$$
In particular,
$\,\soc\Bigl(\td^{\lam}\!\!\downarrow_{\H(D_{n-1})}\Bigr)$ is
multiplicity free.
\end{thm}

\noindent {Proof:}\, This follows directly from Lemma 1.3,
Lemma 1.4, Lemma 2.1, Lemma 2.4 and the fact that
$$\soc\Bigl(\td^{\lam}\!\!\downarrow_{\H(D_{n-1})}\Bigr)=\Bigl\{\soc\Bigl(\td^{\lam}\!\!\downarrow_{\H(B_{n-1})}\Bigr)\Bigr\}\!\!\downarrow_{\H(D_{n-1})}.\eqno(2.5.1)$$

We have to prove (2.5.1). It is clear that $$
\soc\Bigl(\td^{\lam}\!\!\downarrow_{\H(D_{n-1})}\Bigr)\supseteq\Bigl\{\soc\Bigl(\td^{\lam}\!\!\downarrow_{\H(B_{n-1})}\Bigr)\Bigr\}\!\!\downarrow_{\H(D_{n-1})}. $$ It suffices to show that $$
\dim\soc\Bigl(\td^{\lam}\!\!\downarrow_{\H(D_{n-1})}\Bigr)\leq\dim\Bigl\{\soc\Bigl(\td^{\lam}\!\!\downarrow_{\H(B_{n-1})}\Bigr)\Bigr\}\!\!\downarrow_{\H(D_{n-1})}. $$ Since every simple $\H(D_{n-1})$-module occurs as a direct summand of 
$\td^{\mu}\!\!\downarrow_{\H(D_{n-1})}$ for some $\mu\in\mathcal{P}_{n-1}$, we divide the proof into two cases:
\smallskip

\noindent {\it Case 1.} Let $\mu\in\mathcal{P}_{n-1}$ be such that
$\mu\neq\widehat{\mu}$ and
$$\Hom_{\H(D_{n-1})}\Bigl(\td^{\mu}\!\!\downarrow_{\H(D_{n-1})},\td^{\lam}\!\!\downarrow_{\H(D_{n-1})}\Bigr)\neq 0. $$
Then by Frobenius Reciprocity (\cite[(11.13)]{CR}),
$$\aligned
&\quad\,
\Hom_{\H(D_{n-1})}\biggl(\td^{\mu}\!\!\downarrow_{\H(D_{n-1})},\Bigl\{\soc\Bigl(\td^{\lam}\!\!\downarrow_{\H(B_{n-1})}\Bigr)\Bigr\}\!\!\downarrow_{\H(D_{n-1})}\biggr)\\
&\cong\Hom_{\H(B_{n-1})}\biggl(\bigl(\td^{\mu}\!\!\downarrow_{\H(D_{n-1})}\bigr)\!\!\uparrow^{\H(B_{n-1})},
\soc\Bigl(\td^{\lam}\!\!\downarrow_{\H(B_{n-1})}\Bigr)\biggr)\\
&\cong\Hom_{\H(B_{n-1})}\biggl(\td^{\mu}\oplus\bigl(\td^{\mu}\bigr)^{\sigma},
\soc\Bigl(\td^{\lam}\!\!\downarrow_{\H(B_{n-1})}\Bigr)\biggr)\\
&\cong\Hom_{\H(B_{n-1})}\Bigl(\td^{\mu}\oplus\bigl(\td^{\mu}\bigr)^{\sigma},\td^{\lam}\!\!\downarrow_{\H(B_{n-1})}\Bigr)\\
&\cong
\Hom_{\H(B_{n-1})}\Bigl(\bigl(\td^{\mu}\!\!\downarrow_{\H(D_{n-1})}\bigr)\!\!\uparrow^{\H(B_{n-1})},
\td^{\lam}\!\!\downarrow_{\H(B_{n-1})}\Bigr)\\
&\cong\Hom_{\H(D_{n-1})}\Bigl(\td^{\mu}\!\!\downarrow_{\H(D_{n-1})},\td^{\lam}\!\!\downarrow_{\H(D_{n-1})}\Bigr),
\endaligned$$
as required.
\smallskip

\noindent {\it Case 2.} Let $\mu\in\mathcal{P}_{n-1}$ be such that
$\mu=\widehat{\mu}$ and
$$\Hom_{\H(D_{n-1})}\Bigl(D_{+}^{\mu},\td^{\lam}\!\!\downarrow_{\H(D_{n-1})}\Bigr)\neq 0. $$
Then by Frobenius Reciprocity (\cite[(11.13)]{CR}),
$$\aligned &\quad\,
\Hom_{\H(D_{n-1})}\biggl(D_{+}^{\mu},\Bigl\{
\soc\Bigl(\td^{\lam}\!\!\downarrow_{\H(B_{n-1})}\Bigr)\Bigr\}\!\!\downarrow_{\H(D_{n-1})}\biggr)\\
&\cong\Hom_{\H(B_{n-1})}\biggl(D_{+}^{\mu}\!\!\uparrow^{\H(B_{n-1})},\soc\Bigl(\td^{\lam}\!\!\downarrow_{\H(B_{n-1})}\Bigr)\biggr)\\
&\cong\Hom_{\H(B_{n-1})}\biggl(\td^{\mu},\soc\Bigl(\td^{\lam}\!\!\downarrow_{\H(B_{n-1})}\Bigr)\biggr)\\
&\cong\Hom_{\H(B_{n-1})}\Bigl(\td^{\mu},\td^{\lam}\!\!\downarrow_{\H(B_{n-1})}\Bigr)\\
&\cong
\Hom_{\H(B_{n-1})}\Bigl(D_{+}^{\mu}\!\!\uparrow^{\H(B_{n-1})},\td^{\lam}\!\!\downarrow_{\H(B_{n-1})}\Bigr)\\
&\cong\Hom_{\H(D_{n-1})}\Bigl(D_{+}^{\mu},\td^{\lam}\!\!\downarrow_{\H(D_{n-1})}\Bigr),\endaligned
$$
as required. This completes the proof of (2.5.1).\hfill\qed
\medskip

Let $\lam\in\mathcal{P}_n$. Suppose that $\lam\neq\widehat{\lam}$
and $\lam$ is not almost symmetric. We claim that for any two $l$-good
nodes $B, C$ of $\lam$,
$\lam\setminus{B}\neq\widehat{(\lam\setminus{C})}$. In fact, it is
enough to show that if $B\neq C$, then
$\lam\setminus{B}\neq\widehat{(\lam\setminus{C})}$. Otherwise, if
$B, C$ both lie in the same component of $\lam$, say,
$\lam^{(1)}$, then we have that
$\lam^{(1)}\setminus{B}=\lam^{(2)}=\lam^{(1)}\setminus{C}$, which
is impossible; while if $B, C$ lie in different components of $\lam$,
say, $B\in\lam^{(1)}, C\in\lam^{(2)}$, then
$\lam\setminus{B}=\widehat{(\lam\setminus{C})}$ implies that
$\lam^{(1)}=\lam^{(2)}$, which is again impossible. This proves
our claim. It follows from this and Lemma 1.3 and Lemma 1.4
and Lemma 2.1 and (2.5.1) that

\begin{thm} Let $\lam\in\mathcal{P}_n$. Suppose that
$\lam\neq\widehat{\lam}$ and $\lam$ is not almost symmetric. Then
$$ \soc\Bigl(\td^{\lam}\!\!\downarrow_{\H(D_{n-1})}\Bigr)\cong
\bigoplus_{\substack{C\in [\lam]\\ \text{$C$ is $l$-good}}}
\td^{\lam\setminus{C}}\!\!\downarrow_{\H(D_{n-1})}.
$$
In particular,
$\,\soc\Bigl(\td^{\lam}\!\!\downarrow_{\H(D_{n-1})}\Bigr)$ is
multiplicity free.
\end{thm}

Let $\lam\in\mathcal{P}_n$. Now suppose that $\lam=\widehat\lam$.
It remains to describe the decompositions of the socle of
$D_+^{\lam}\!\!\downarrow_{\H(D_{n-1})}$ and of
$D_-^{\lam}\!\!\downarrow_{\H(D_{n-1})}$ into irreducible
$\H(D_{n-1})$-modules.

\begin{thm} Let $\lam\in\mathcal{P}_n$ be such that
$\lam=\widehat\lam$. Then there is a $\H(D_{n-1})$-module
isomorphism $$
\soc\Bigl({D_+^{\lam}}\!\!\downarrow_{\H(D_{n-1})}\Bigr)\cong\soc\Bigl({D_-^{\lam}}\!\!\downarrow_{\H(D_{n-1})}\Bigr).
$$
\end{thm}

\noindent {Proof:}\, By assumption $n$ is even. Hence $n-1$ is
odd. In particular, for any bipartition $\mu\in\mathcal{P}_{n-1}$,
$\mu\neq\widehat{\mu}$. By \cite[Corollary 2.4]{Hu3},
$\td^{\mu}\cong\Bigl(\td^{\mu}\Bigr)^{\tau}$. We have that
$$\aligned
&\quad\,\Hom_{\H(D_{n-1})}\Bigl(\td^{\mu}\!\!\downarrow_{\H(D_{n-1})},{D_+^{\lam}}\!\!\downarrow_{\H(D_{n-1})}\Bigr)\\
&\cong
\Hom_{\H(D_{n-1})}\Bigl(\bigl(\td^{\mu}\!\!\downarrow_{\H(D_{n-1})}\bigr)^{\tau},\bigl({D_+^{\lam}}\!\!\downarrow_{\H(D_{n-1})}\bigr)^{\tau}\Bigr)\\
&\cong\Hom_{\H(D_{n-1})}\Bigl(\bigl(\td^{\mu}\bigr)^{\tau}\!\!\downarrow_{\H(D_{n-1})},\bigl({D_+^{\lam}}\bigr)^{\tau}\!\!\downarrow_{\H(D_{n-1})}\Bigr)\\
&\cong\Hom_{\H(D_{n-1})}\Bigl(\td^{\mu}\!\!\downarrow_{\H(D_{n-1})},{D_-^{\lam}}\!\!\downarrow_{\H(D_{n-1})}\Bigr).
\endaligned$$ Now using Lemma 1.3 and Lemma 1.4, the
theorem follows at once.  \hfill\qed\medskip

We define an equivalence relation $\sim$ on ${\mathcal P}_{n-1}$
by $\lam\sim\mu$ if and only if $\mu=\widehat\lam$. Then Lemma
2.1 and (2.5.1) implies the following:

\addtocounter{cor}{7}
\begin{cor} Let $\lam\in\mathcal{P}_n$ be such that
$\lam=\widehat\lam$. Then $$
\soc\Bigl({D_+^{\lam}}\!\!\downarrow_{\H(D_{n-1})}\Bigr)\cong
\soc\Bigl({D_-^{\lam}}\!\!\downarrow_{\H(D_{n-1})}\Bigr)
\cong\bigoplus_{\mu} \td^{\mu}\!\!\downarrow_{\H(D_{n-1})},
$$
where the sum $\mu$ is taken over a fixed set of representatives
of equivalence classes in $\mathcal{P}_{n-1}/{\sim}$ such that
$\mu\rightarrow\lam$. In particular,
$\,\soc\Bigl(D_+^{\lam}\!\!\downarrow_{\H(D_{n-1})}\Bigr)$ and
$\soc\Bigl(D_-^{\lam}\!\!\downarrow_{\H(D_{n-1})}\Bigr)$ are both
multiplicity free.
\end{cor}

\begin{cor} For any irreducible $\H(D_n)$-module $D$,
$\,\soc\Bigl(D\!\!\downarrow_{\H(D_{n-1})}\Bigr)$ is multiplicity
free.
\end{cor}

\noindent Now Theorem 2.5, Theorem 2.6 and Corollary 2.8
completely determine the decomposition of
$\soc\Bigl(D\!\!\downarrow_{\H(D_{n-1})}\Bigr)$ into irreducible
$\H(D_{n-1})$-modules for every irreducible $\H(D_n)$-module $D$.

\addtocounter{example}{6}
\begin{example} Suppose $\H(D_5)$ is semi-simple. With
the notations in Example 2.3, we have that $$ \begin{aligned}
\ts^{\lam}\!\!\downarrow_{\H(D_4)}&\cong\ts_+^{((1^2),(1^2))}\bigoplus
\ts_-^{((1^2),(1^2))}\bigoplus\ts^{((2),(1^2))}\!\!\downarrow_{\H(D_4)}\\
&\qquad\qquad
\bigoplus\ts^{((2,1),(1))}\!\!\downarrow_{\H(D_4)},\\
\ts^{\mu}\!\!\downarrow_{\H(D_4)}&\cong\ts^{((1),(1^3))}\!\!\downarrow_{\H(D_4)}\bigoplus
\ts^{((2),(1^2))}\!\!\downarrow_{\H(D_4)}.\end{aligned}$$
\end{example}

\begin{example} Suppose  $\H(D_6)$ is semi-simple. Let
$n=6$, $\nu=((2,1), (2,1))$. Then we have that $$
\ts_+^{\nu}\!\!\downarrow_{\H(D_5)}\cong\ts_-^{\nu}\!\!\downarrow_{\H(D_5)}\cong
\ts^{((2,1),(2))}\!\!\downarrow_{\H(D_5)}\bigoplus
\ts^{((2,1),(1^2))}\!\!\downarrow_{\H(D_5)}.
$$
\end{example}

\begin{example} Suppose that $q=1$. Then, as long as
$\ch K\neq 2$, Theorem 2.5, Theorem 2.6 and Corollary 2.8
completely determine the decomposition of
$\soc\Bigl(D\!\!\downarrow_{KW(D_{n-1})}\Bigr)$ into irreducible
$KW(D_{n-1})$-modules for every modular irreducible
$KW(D_n)$-module $D$.
\end{example}

\section{The case where
$\prod_{i=1}^{n-1}\bigl(1+q^i\bigr)=0$ and $2\cdot 1_K\neq 0$}

Throughout this section, we assume that
$\prod_{i=1}^{n-1}\bigl(1+q^i\bigr)=0$,  $2\cdot 1_K\neq 0$ and
$K$ is a field such that $\H(D_n)$ is split over $K$. It follows
that $q$ is a primitive $2\l$-th root of unity for some integer
$1\leq\l<n$. In this case $\mathcal{P}_n$ can be identified with
the set of Kleshchev bipartitions of $n$  with respect to
$(\sqrt[2\l]{1};1,-1)$ (see \cite{A}).

Let $v$ be an indeterminate over $\mathbb{Q}$. Let $\mathfrak{h}$
be a $(2\l+1)$-dimensional vector space over $\mathbb{Q}$ with
basis $\{h_0,h_1,\cdots,$ $h_{2\l-1},d\}$. Denote by
$\{\Lambda_0,\Lambda_1,\cdots,\Lambda_{2\l-1},\delta\}$ the dual
basis of $\mathfrak{h}^*$, and we set
$\alpha_i=2\Lambda_i-\Lambda_{i-1}-\Lambda_{i+1}+\delta_{i,0}\delta$
for $i\in\Z/{2\l\Z}$. Assume that the $2\l\times 2\l$ matrix
$(\langle\alpha_i,h_j\rangle)$ is the generalized Cartan matrix
associated to $\ksl_{2\l}$. Let $U_v(\ksl_{2\l})$ be the quantum
affine algebra corresponding to $\ksl_{2\l}$ (see \cite[\S2]{Hu2}
for its definition). Let $\Lambda=\Lambda_0+\Lambda_{\l}$. Let
$\mathcal{F}(\Lambda):=\oplus_{n\geq
0}\oplus_{\lam\in\widetilde{\mathcal{P}}_n}\mathbb{Q}(v)\lam$,
{\it a level $2$ Fock space}, on which $U_v(\ksl_{2\l})$ acts. The
submodule $L(\Lambda)$ generated by the empty bipartition $\underline{\emptyset}$ is the
irreducible integrable highest weight module with highest weight
$\Lambda$. By a well-known result of S. Ariki (\cite{A2}), the dual of the
Grothendieck group $K(\oplus_{k\geq 0}\H(B_k))$ can be made into a
$U_v(\ksl_{2\l})$-module. Ariki introduced the functors of $i$-restriction
and $i$-induction, which plays the role of Chevalley generators $e_i, f_i$.
It is a remarkable fact (see \cite{MM},
\cite[(2.11)]{AM}) that the crystal graph of $L(\Lambda)$ can be realized as 
the Kleshchev's good lattice with respect to $(\sqrt[2\l]{1};1,-1)$
if one uses the embedding $L(\Lambda)\subset\mathcal{F}(\Lambda)$. 
The operators of removing and adding $2\l$-good nodes play the role of the Kashiwara operators 
$\tilde{e}_i,\tilde{f}_i$. On the other hand, Grojnowski (see \cite[Theorem 14.2, 14.3]{G})
gives another realization of the crystal graph of $L(\Lambda)$ on the set of all the simple modules of 
$\H(B_k)$ for all $k\geq 0$, where the functors of taking socle (resp., taking cosocle) of the $i$-restriction
(resp. of the $i$-induction) of simple modules play the role of the Kashiwara operators 
$\tilde{e}_i,\tilde{f}_i$. These two crystal structures are isomorphic to each other. 
By the definition of the second realization, there are the following results (see \cite{Ma}).

\begin{lem}\text{(\cite{G},\cite{GV})} There is an isomorphism $\pi$ between the above two realizations of the crystal structure, such that, if we write $\pi=\oplus_{k\geq 0}\pi_k$, where
$\pi_k$ is a permutation defined on $\mathcal{P}_k$, then for each
bipartition $\lam\in\mathcal{P}_n$, $$
\soc\Bigl(\td^{\pi_n(\lam)}\!\!\downarrow_{\H(B_{n-1})}\Bigr)\cong\bigoplus_{\mu\rightarrow\lam}
\td^{\pi_{n-1}(\mu)}, $$ where $\mu\rightarrow\lam$ means
that $\mu$ is a bipartition of $n-1$ such that the Young diagram
$[\mu]$ is obtained from the Young diagram $[\lam]$ by removing a
$2\l$-good node. In particular, for any bipartition
$\lam\in\mathcal{P}_n$, the socle of
$\td^{\lam}\!\!\downarrow_{\H(B_{n-1})}$ is a direct sum of pairwise
non-isomorphic irreducible $\H(B_{n-1})$-modules, i.e., it is
multiplicity free.
\end{lem}

We remark that both Ariki's and Grojnowski's results are stated in the context of general cyclotomic Hecke algebras of type $G(r,1,n)$, though we only use the special type $B$ case in this paper. It has been conjectured for a long time that the above $\pi$ can be chosen as identity (compare \cite{Kl}, \cite{B}). In a recent preprint \cite{A3}, Ariki proved this conjecture. As a consequence, we have the following.

\begin{lem} [\rm \cite{A3}] With the above
notations, $\pi_n$ can be chosen as the identity map for any $n$.
\end{lem}

To describe the decomposition into irreducible modules of the
socle of the module $D\!\!\downarrow_{\H(D_{n-1})}$ for each
irreducible $\H(D_n)$-module $D$ in this case, we need a better
understanding of the set of fixed-points under the involution
$\bH$.

\begin{lem} Let $\lam\in\mathcal{P}_n$ be a Kleshchev
bipartition with respect to $(\sqrt[2\l]{1};1,-1)$. Suppose that
$\lam=\bH(\lam)$. Then $n$ is even. Moreover, if
$\emptyset\stackrel{i_1}{\longrightarrow}\cdot\overset{i_2}{\longrightarrow}\cdot
\cdots\stackrel{i_n}{\longrightarrow}\lam $ is a path in
Kleshchev's good lattice with respect to $(\sqrt[2\l]{1};1,-1)$, then for any
integer $k$, we have $$ \#\bigl\{j\bigm| \text{$i_j\equiv
k\!\!\!\mod{2\l}$}\bigr\}=\#\bigl\{j\bigm|\text{$i_j\equiv{k+\l}\!\!\!\mod{2\l}$}\bigr\}.
$$
\end{lem}

\noindent {Proof:} \, Let $$
\emptyset\stackrel{i_1}{\longrightarrow}\cdot\stackrel{i_2}{\longrightarrow}\cdot
\cdots\stackrel{i_n}{\longrightarrow}\lam $$ be a path in
Kleshchev's good lattice with respect to $(\sqrt[2\l]{1};1,-1)$. By assumption
and \cite[3.3]{Hu2}, $$
\emptyset\stackrel{i_1+\l}{\longrightarrow}\cdot\stackrel{i_2+\l}{\longrightarrow}\cdot
\cdots\stackrel{i_n+\l}{\longrightarrow}\lam
$$
is also a path in Kleshchev's good lattice with respect to
$(\sqrt[2\l]{1};1,-1)$.

It follows that there is an automorphism defined on the set of
nodes of $\lam$, say $\psi$, such that
$\res\gamma=\res\psi(\gamma)+\l$ in $\Z/{2\l\Z}$. We claim that
this is enough for deducing our theorem by using induction on $n$.

In fact, we first pick a node $A_1\in\lam$ and define
$B_1:=\psi(A_1)$. As $\res A_1=\res B_1+\l\neq\res B_1$, it is
clear that $A_1\neq B_1$. Now if $\psi(B_1)=A_1$, then we can
remove $A_1, B_1$ and use induction on $n$; otherwise, we define
$A_2:=\psi(B_1)$, which is different from $A_1$. Since $\psi$ is
bijective, it follows that $B_2:=\psi(A_2)$ is different from
$B_1=\psi(A_1)$. Repeating this procedure, and noting that $\lam$
has only finitely many nodes, one would get $2k$ pairwise
different nodes of $\lam$, $A_1, A_2,\cdots, A_k$ and
$B_1,B_2,\cdots, B_k$ such that
$B_i=\psi(A_i),\,\,\forall\,\,1\leq i\leq k$ and
$A_{j}=\psi(B_{j-1}),\,\,\forall\,2\leq j\leq k$ and
$\psi(B_{k})=A_{1}$. Then we can remove these $2k$ nodes and use
induction hypothesis on $n$. This completes the proof of our
theorem. \hfill\qed

\addtocounter{defs}{3}
\begin{defs} Let $\lam\in\mathcal{P}_n$. Suppose that
$\lam$ has a $2\l$-good node $A$ such that $$
\lam\setminus{A}=\bH(\lam\setminus{A}).
$$
Then we say that $\lam$ is almost $\l$-symmetric.
\end{defs}

\addtocounter{example}{4}
\begin{example} Suppose that $n=5$ and $\l=2$. Let
$\lam=((1), (2^2)),\, \mu=((2), (1^3))$ are two Kleshchev
bipartitions of $5$ with respect to $(\sqrt[4]{1};1,-1)$. Let $A$ be the
node which is in the second row and the second column of the
second component of $\lam$. Then $\lam$ is almost $2$-symmetric
with $\lam\setminus{A}=\bH(\lam\setminus{A})$, but $\mu$ is not
almost $2$-symmetric.
\end{example}

\addtocounter{lem}{2}
\begin{lem} Let $\lam\in\mathcal{P}_n$. Suppose that
$\lam$ is almost $\l$-symmetric with
$\,\lam\setminus{A}=\bH(\lam\setminus{A})$ for some $2\l$-good node $A$
of $\,\lam$. Then for any pairs of $2\l$-good nodes $B, C$ of $\,\lam$
satisfying $C\neq A$, we have that $$
\lam\setminus{B}\neq\bH(\lam\setminus{C}).
$$
In particular, for any $2\l$-good nodes $C$ of $\,\lam$ satisfying
$C\neq A$, we have that $$
\lam\setminus{C}\neq\bH(\lam\setminus{C}). $$
\end{lem}

\noindent {Proof:} \, Suppose that
$\lam\setminus{B}=\bH(\lam\setminus{C})$. Since $\bH$ is an
involution, $\,\lam\setminus{A}=\bH(\lam\setminus{A})$ and $C\neq
A$, it follows that $B\neq A$. Note that for $2\l$-good nodes $A, C$
(resp. $A, B$), $A\neq C$ (resp. $A\neq B$) implies that their
residues $\res A, \res C$ (resp. $\res A, \res B$) are different.
Write $\res A=i$, $\res B=j$, $\res C=k$. Then $j\neq i\neq k$.

For each bipartition $\mu$ and each $s\in\Z/{2\l\Z}$, we denote by
$N_s^{\mu}$ the number of $s$-nodes in the Young diagram $[\mu]$
of $\mu$. If $j\neq k+\l$, then by Lemma 1.4(2),
$$N_{k+\l}^{\lam\setminus{B}}=N_{k+\l}^{\lam}\geq
N_{k+\l}^{\lam\setminus{A}}=N_{k+\l}^{\bH(\lam\setminus{A})}=N_{k}^{\lam\setminus{A}}=N_{k}^{\lam\setminus{C}}+1=
N_{k+\l}^{\bH(\lam\setminus{C})}+1, $$ a contradiction. If $j=
k+\l$, then by Lemma 1.4(2),$$
N_{i}^{\lam\setminus{B}}=N_{i}^{\lam}=
N_{i}^{\lam\setminus{A}}+1=N_{i}^{\bH(\lam\setminus{A})}+1=N_{i+\l}^{\lam\setminus{A}}+1=N_{i+\l}^{\lam\setminus{C}}+1=
N_{i}^{\bH(\lam\setminus{C})}+1,
$$
a contradiction. This proves our lemma. \hfill\qed

\addtocounter{thm}{6}
\begin{thm} Let $\lam\in\mathcal{P}_n$. Suppose that
$\lam$ is almost $\l$-symmetric with
$\lam\setminus{A}=\bH(\lam\setminus{A})$ for some $2\l$-good node $A$ of
$\lam$. Then $\lam\neq\bH(\lam)$ and $$
\soc\Bigl(D^{\lam}\!\!\downarrow_{\H(D_{n-1})}\Bigr)\cong
D_{+}^{\lam\setminus{A}}\bigoplus
D_{-}^{\lam\setminus{A}}\bigoplus
\bigoplus_{\substack{C\in [\lam],\,\, C\neq A\\ \text{$C$
is $2\l$-good}}}
\td^{\lam\setminus{C}}\!\!\downarrow_{\H(D_{n-1})}.
$$
In particular, $\,\td^{\lam}\!\!\downarrow_{\H(D_{n-1})}$ is
multiplicity free.
\end{thm}

\noindent {Proof:}\, This follows directly from Lemma 1.3,
(2.5.1), Lemma 3.1, Lemma 3.2, Lemma 3.3 and Lemma 3.6. \hfill\qed

\begin{thm} Let $\lam\in\mathcal{P}_n$. Suppose that
$\lam\neq\bH(\lam)$ and $\lam$ is not almost $\l$-symmetric. Then
for any two $2\l$-good nodes $B, C$ of $\lam$,
$\lam\setminus{B}\neq\bH(\lam\setminus{C})$.
\end{thm}

\noindent {Proof:} \, In fact, it suffices to show that, if $B\neq
C$, then $\lam\setminus{B}\neq\bH(\lam\setminus{C})$.

Otherwise, suppose that $B\neq C$ and
$\lam\setminus{B}=\bH(\lam\setminus{C})$. We call two nodes
$\gamma, \gamma'$ $\l$-conjugate, if $\res\gamma=\res\gamma'+\l$
in $\Z/{2\l\Z}$. We deduce that $B$ is not $\l$-conjugate to $C$
(otherwise it would follow that $\lam=\bH(\lam)$). Therefore, we
have that $\res B\neq\res C$ (as $B\neq C$ are both $2\l$-good nodes),
and $\res B\neq\res C+\l$. Now the condition that
$\lam\setminus{B}=\bH(\lam\setminus{C})$ implies that there is a
bijection, say $\varphi$, from the set of the nodes of
$\lam\setminus{B}$ onto the set of the nodes of
$\lam\setminus{C}$, such that for any
$\gamma\in{\lam\setminus{B}}$, $\gamma$ is $\l$-conjugate to
$\varphi(\gamma)$.

We define $C_1=C_0:=C$, $D_1:=\varphi(C_1)$. Since $D_1$ is
$\l$-conjugate to $C_1$ and hence different from $B$, we have that
$D_1\in\lam\setminus{B}$. Hence we can define
$C_2:=\varphi(D_1)\in\lam\setminus{C}$, then $\res C_2=\res C_1$
and hence $C_2\neq B$, then we can still define
$D_2:=\varphi(C_2)$. Since $\varphi$ is bijective, $C_1\neq C_2$
implies that $D_1\neq D_2$. In general, suppose that for integer
$k\geq 2$, $D_i\in\lam\setminus{B},\,\, C_{i}\in\lam\setminus{C}$,
for $1\leq i\leq k$, are already well-defined, such that
$C_{j}=\varphi(D_{j-1}),\,\, D_{j}=\varphi(C_{j}),\,\,
\forall\,2\leq j\leq k$ and $$ C_i\neq C_j,\,\,  D_i\neq D_j,
\,\,\res C_i=\res C_j=\res C,\,\, \res D_i=\res D_j=\res C+\l, $$
for any $1\leq i\neq j\leq k$, then we define
$C_{k+1}:=\varphi(D_k)\in\lam\setminus{C}$. As $D_1,\cdots,D_k$
are pairwise different and $\varphi$ is bijective, we get that
$C_1,\cdots,C_{k+1}$ are also pairwise different. Moreover, $\res
C_{k+1}=\res C$, and hence $C_{k+1}\neq B$ and we can still define
$D_{k+1}:=\varphi(C_{k+1})\in\lam\setminus{B}$. It is clear that
$\res D_{k+1}=\res C+\l$, and since $C_1,\cdots,C_{k+1}$ are
pairwise different and $\varphi$ is bijective, we get that
$D_1,\cdots,D_{k+1}$ are also pairwise different. As a
consequence, we get infinitely many pairwise different nodes
$C_1,C_2,\cdots$ in $\lam$, which is impossible. This proves the
theorem.  \hfill\qed

It follows from this theorem,  Lemma 1.3, Lemma 3.1, Lemma 3.2 and (2.5.1) that

\begin{thm} Let $\lam\in\mathcal{P}_n$. Suppose that
$\lam\neq\bH(\lam)$ and $\lam$ is not almost $\l$-symmetric. Then
$$
\soc\Bigl(\td^{\lam}\!\!\downarrow_{\H(D_{n-1})}\Bigr)\cong
\bigoplus_{\substack{C\in [\lam]\\ \text{$C$ is $2\l$-good}}}
\td^{\lam\setminus{C}}\!\!\downarrow_{\H(D_{n-1})}.
$$
In particular,
$\,\soc\Bigl(\td^{\lam}\!\!\downarrow_{\H(D_{n-1})}\Bigr)$ is
multiplicity free.
\end{thm}

Now let $\lam\in\mathcal{P}_n$ be such that $\lam=\bH(\lam)$. It remains to
describe the decompositions into irreducible $\H(D_{n-1})$-modules
of $\soc\Bigl(D_+^{\lam}\!\!\downarrow_{\H(D_{n-1})}\Bigr)$ and
$\soc\Bigl(D_-^{\lam}\!\!\downarrow_{\H(D_{n-1})}\Bigr)$.

\begin{thm} Let $\lam\in\mathcal{P}_n$. Suppose that
$\lam=\bH(\lam)$. Then there is a $\H(D_{n-1})$-module isomorphism
$$ \soc\Bigl(D_+^{\lam}\!\!\downarrow_{\H(D_{n-1})}\Bigr)\cong
\soc\Bigl(D_-^{\lam}\!\!\downarrow_{\H(D_{n-1})}\Bigr).
$$
\end{thm}

\noindent {Proof:}\, By assumption and Lemma 3.3, we know that $n$ is even.
Hence $n-1$ is odd. In particular, for any bipartition
$\mu\in\mathcal{P}_{n-1}$, $\mu\neq\bH(\mu)$. By \cite[Corollary
2.4]{Hu3}, $\td^{\mu}\cong\Bigl(\td^{\mu}\Bigr)^{\tau}$. 
Now using the same argument as in the proof of Theorem
2.7, we prove the theorem. \hfill\qed

\addtocounter{cor}{10}
\begin{cor} Let $\lam\in\mathcal{P}_n$ be such that
$\lam=\bH(\lam)$. In particular, $n$ is even. Then $$
\soc\Bigl(D_+^{\lam}\!\!\downarrow_{\H(D_{n-1})}\Bigr)\cong
\soc\Bigl(D_-^{\lam}\!\!\downarrow_{\H(D_{n-1})}\Bigr)
\cong\bigoplus_{\mu} \td^{\mu}\!\!\downarrow_{\H(D_{n-1})},
$$
where the subscript $\mu$ is taken over a fixed set of
representatives of equivalence classes in
$\mathcal{P}_{n-1}/{\approx}$ such that
$\mu\rightarrow\lam$. In particular,
$\,\soc\Bigl(D_+^{\lam}\!\!\downarrow_{\H(D_{n-1})}\Bigr)$ and
$\,\soc\Bigl(D_-^{\lam}\!\!\downarrow_{\H(D_{n-1})}\Bigr)$ are
both multiplicity free.
\end{cor}

\noindent {Proof:} \, Since $\lam=\bH(\lam)$, by Lemma 3.3, it
follows that for any $2\l$-good node $C$ of $\lam$, $\lam\setminus{C}\neq\bH(\lam\setminus{C})$.
Now using Lemma 3.1, lemma 3.2, (2.5.1) and Theorem 3.10, we prove
the corollary. \hfill\qed

\begin{cor} For any irreducible $\H(D_n)$-module
$D$, the socle of $D\!\!\downarrow_{\H(D_{n-1})}$ is always multiplicity
free.
\end{cor}

Now Theorem 3.7, Theorem 3.9 and Corollary 3.11
completely determine the decomposition of
$\soc\Bigl(D\!\!\downarrow_{\H(D_{n-1})}\Bigr)$ into irreducible
$\H(D_{n-1})$-modules for every irreducible $\H(D_n)$-module $D$.

\section{Appendix}

In this appendix, we give a proof to show that the involution
$\bH$ (and hence the main result of \cite{Hu2}) is independent of
the base field $K$ as long as $\ch K\neq 2$ and $\H(D_n)$ is split
over $K$. The proof is essentially due to Professor S. Ariki.

Throughout this section, we assume that $\ch K\neq 2$ and
$\H(D_n)$ is split over $K$. By Lemma 1.4(1), it suffices to
consider the case where $q$ is a primitive $2\l$-th root of unity
in $K$ for some integer $1\leq\l<n$. To emphasize the base field
$K$, we denote by $\H_K(D_n)$ the Hecke algebra of type $D_n$ over
$K$, and by $\ts_K^{\lam}$ (resp. $\td_K^{\lam}$) the
corresponding $\H_K(D_n)$-modules. Note that by \cite{A}, the set
$\bigl\{\lam\in\widetilde{\mathcal{P}}_n\bigm|\td_K^{\lam}\neq
0\bigr\}$ depends only on $\l$, but not on the choice of the base
field $K$. So we denote it by $\mathcal{P}_n$ as before. There is
an involution $\bH_K$ defined on the set $\mathcal{P}_n$ such
that, for each $\lam\in\mathcal{P}_n$,
$\Bigl(\td_K^{\lam}\Bigr)^{\sigma}\cong\td_K^{\bH_K(\lam)}$. We
have that

\begin{thm} Suppose that $\ch K=p\neq 2$ and
$\H_K(D_n)$ is split over $K$. Then $\bH_{K}=\bH_{\mathbb{C}}$.
\end{thm}

\noindent {Proof:} \, By \cite{DJMu}, we know that each
irreducible $\H(B_n)$-module remains irreducible under field
extension. By definition of the automorphism $\sigma$ and the
involution $\bH$, it is easy to see that if
$\bH_{F}=\bH_{\mathbb{C}}$ for some splitting field $F$ of
$\H(D_n)$ with $\ch F=p\neq 2$, then for any splitting field $K'$
of $\H(D_n)$ with $\ch K'=p\neq 2$, we have that
$\bH_{K'}=\bH_{\mathbb{C}}$. In particular, for any characteristic
$0$ splitting field $E$ of $\H(D_n)$, we have that
$\bH_{E}=\bH_{\mathbb{C}}$.

Therefore, it suffices to consider the characteristic $p>2$ case.
To ensure the existence of a primitive $2\l$-th root of unity, we
further assume that $(p,\l)=1$. Let $q\in\overline{\mathbb{F}_p}$
(resp. $q_0\in\mathbb{C}$) be a primitive $2\l$-th root of unity,
where $\overline{\mathbb{F}_p}$ is the algebraic closure of the
finite field $\mathbb{F}_p$. Let $X$ be an indeterminate over
$\mathbb Z$. For each polynomial $f\in\mathbb{Z}[X]$, let
$\overline{f}$ be its canonical image in $\mathbb{F}_p[X]$. For
each $m\in\mathbb N$, let $\Phi_m(X)$ be the $m$-th cyclotomic
polynomial over $\mathbb Z$. Then $X^m-1=\prod_{d\mid
m}\Phi_d(X)$. It follows that $\overline{\Phi_{2\l}}(q)=0$. Hence
the map which sends $q_0$ to $q$ extends naturally to a surjective
ring homomorphism from $\mathbb{Z}[q_0]$ onto $\mathbb{F}_p[q]$.
Let $\Phi_{2\l,p}(X)$ be a monic polynomial in $\mathbb{Z}[X]$
such that $\overline{\Phi_{2\l,p}}(X)$ is the minimal polynomial
of $q$ over $\mathbb{F}_p$.

Recall that every finite dimensional algebra becomes split after a
finite field extension. Therefore, there exist some algebraic
integers $\alpha_1,\cdots,\alpha_s\in\mathbb{C}$, some elements
$\overline{\alpha_1},\cdots,\overline{\alpha_s}\in
\overline{\mathbb{F}_p}$, and some monic polynomials, say
$f_i(X_i)\in\mathbb{Z}[q_0,\alpha_1,\cdots,$ $\alpha_{i-1}][X_i]$,
$\overline{f_i}(X_i)\in\mathbb{F}_p[q,\overline{\alpha_1},\cdots,
\overline{\alpha_{i-1}}][X_i]$, where $X_i$ is an indeterminate
over $\mathbb{Z}[q_0,\alpha_1,\cdots,\alpha_{i-1}]$ (resp. over
$\mathbb{F}_p[q,\overline{\alpha_1},\cdots,
\overline{\alpha_{i-1}}]$), $i=1,\cdots,s$, such that
\smallskip

1) $\overline{f_i}$ (respectively, $f_i$) is irreducible over
$\mathbb{F}_p[q,\overline{\alpha_1},
\cdots,\overline{\alpha_{i-1}}]$ (respectively, over
$\mathbb{Z}[q_0,\alpha_1,\cdots,$ $\alpha_{i-1}]$),

2) $f_i(\alpha_i)=0$, $\overline{f_i}(\overline{\alpha_i})=0$, and
$\H_{{\mathbb{F}_p}[q,\overline{\alpha_1},\cdots,\overline{\alpha_s}]}(D_n)$
is split over the field
${\mathbb{F}_p}[q,\overline{\alpha_1},\cdots,\overline{\alpha_s}]$.
\smallskip

Note that (see \cite[Chapter IV, \S1, Theorem 4]{L})
$\mathbb{Z}[q_0]$ is the integral closure of $\mathbb Z$ in
$\mathbb{Q}[q_0]$. Let $R'$ be the integral closure of $\mathbb
Z[q_0]$ in $\mathbb{Q}[q_0,\alpha_1,\cdots,\alpha_{s}]$. By
\cite[Chapter 5, Exercise 2]{AtM}, the natural surjective morphism
from $\mathbb Z[q_0,\alpha_1,\cdots,\alpha_{s}]$ onto
${\mathbb{F}_p}[q,\overline{\alpha_1},\cdots,\overline{\alpha_s}]\subseteq
\overline{\mathbb{F}_p}$ can be extended to a ring homomorphism
from $R'$ to the field $\overline{\mathbb{F}_p}$. We denote it by
$\pi$. Then (by \cite[Chapter I, \S4]{S}) $R'$ is a Dedekind
domain and the field of fractions of $R'$ is
$\mathbb{Q}[q_0,\alpha_1,\cdots,\alpha_{s}]$.

Similarly, let $E$ be a finite extension of
$\mathbb{Q}[q_0,\alpha_1,\cdots,\alpha_{s}]$ such that $\H_E(D_n)$
is split over the field $E$. Let $R$ be the integral closure of
$R'$ in $E$. By \cite[Chapter I, \S4]{S}, $R$ is a Dedekind domain
and the field of fractions of $R$ is $E$. By \cite[Chapter 5,
Exercise 2]{AtM}, the homomorphism $\pi$ can be extended to a ring
homomorphism from $R$ to the field $\overline{\mathbb{F}_p}$. It
follows that the ideal of $R$ generated by $p$ and
$\Phi_{2\l,p}(q_0)$ should be a proper ideal. Let $\mathfrak{m}$
be the kernel of the homomorphism, which is a maximal ideal of $R$
containing $p$ and $\Phi_{2\l,p}(q_0)$. Let
$\mathcal{O}:=R_{\mathfrak{m}}$,
$F:=R_{\mathfrak{m}}/\mathfrak{m}R_{\mathfrak{m}}$. It is clear
that
${\mathbb{F}_p}[q,\overline{\alpha_1},\cdots,\overline{\alpha_s}]\subseteq
F$. Therefore we get a $p$-modular system $(\mathcal{O}, E, F)$,
where $E$ (resp. $F$) is a field of characteristic $0$ (resp.
characteristic $p$) such that $\H_E(D_n)$ (resp. $\H_F(D_n)$) is
split over $E$ (resp. over $F$), and $q_0\in\mathcal{O}\subset E$
is a primitive $2\l$-th root of unity in $E$ which is in the
pre-image of $q$. By results of \cite{GR} and \cite[(2.3)]{Ge},
the decomposition map from the Grothendieck group $K_0(\H_E(D_n))$
to the Grothendieck group $K_0(\H_F(D_n))$ is well-defined.

Let $\lam\in\mathcal{P}_n$. Recall that there is a well-defined
bilinear form $\langle,\rangle_{\mathcal{O}}$ over
$\ts_{\mathcal{O}}^{\lam}$. Let
$M_{\mathcal{O}}:=\rad\langle,\rangle_{\mathcal{O}}$,
$M:=\rad\langle,\rangle_{E}$. Then $M\cong
M_{\mathcal{O}}\otimes_{\mathcal{O}}E$ is the unique maximal
$\H_E(B_n)$-submodule of $\ts_E^{\lam}$, and $\td_E^{\lam}\cong
 \ts_E^{\lam}/M$. Let $0\neq a\in\mathcal{O}$, $z\in M_{\mathcal{O}}$.
 If $ax=z$ for some $x\in\ts_{\mathcal{O}}^{\lam}$, then $a\langle x,y\rangle=\langle ax,y\rangle=0$
for any $y\in\ts_{\mathcal{O}}^{\lam}$. Since $\mathcal{O}$ is an
integral domain, it follows that $\langle x,y\rangle=0$ for any
$y\in\ts_{\mathcal{O}}^{\lam}$, hence $x\in M_{\mathcal{O}}$. This
shows that $M_{\mathcal{O}}$ is a
 pure $\mathcal{O}$-submodule of $\ts_{\mathcal{O}}^{\lam}$. Since
 $\mathcal{O}$ is a principal integral domain,
it follows\footnote{Recall that over a principal integral domain,
 any pure submodule of a finitely generated module must be its direct
 summand.} that $M_{\mathcal{O}}$ is a
 direct $\mathcal{O}$-summand of $\ts_{\mathcal{O}}^{\lam}$. In
 particular,
 $\td_{\mathcal{O}}^{\lam}:=\ts_{\mathcal{O}}^{\lam}/M_{\mathcal{O}}$
 and hence $\Bigl(\td_{\mathcal{O}}^{\lam}\Bigr)^{\sigma}$
 are both free $\mathcal{O}$-modules.

Note that
$$ \td_E^{\lam}\cong\ts_E^{\lam}/M\cong (\ts_{\mathcal{O}}^{\lam}\otimes_{\mathcal{O}}
E)/(M_{\mathcal{O}}\otimes_{\mathcal{O}}
E)\cong\td_{\mathcal{O}}^{\lam}\otimes_{\mathcal{O}} E,
$$
and hence
$$\td_E^{\bH_{\mathbb{C}}(\lam)}\cong\Bigl(\td_E^{\lam}\Bigr)^{\sigma}\cong\Bigl(\td_{\mathcal{O}}^{\lam}\otimes_{\mathcal{O}}
E\Bigr)^{\sigma} \cong
\Bigl(\td_{\mathcal{O}}^{\lam}\Bigr)^{\sigma}\otimes_{\mathcal{O}}
E,\eqno{(4.2)}$$ which implies that
$\Bigl(\td_{\mathcal{O}}^{\lam}\Bigr)^{\sigma}$ is a full
$\H_{\mathcal{O}}(B_n)$-lattice in
 $\td_{E}^{\bH_{\mathbb{C}}(\lam)}$.
By the fact that
$\td_E^{\bH_{\mathbb{C}}(\lam)}\cong\td_{\mathcal{O}}^{\bH_{\mathbb{C}}(\lam)}\otimes_{\mathcal{O}}
E$, we know that $\td_{\mathcal{O}}^{\bH_{\mathbb{C}}(\lam)}$ is
also a full $\H_{\mathcal{O}}(B_n)$-lattice in
 $\td_{E}^{\bH_{\mathbb{C}}(\lam)}$. Therefore, the module
$\Bigl(\td_{\mathcal{O}}^{\lam}\Bigr)^{\sigma}\otimes_{\mathcal{O}}
F$ has the same set of composition factors as that of
$\td_{\mathcal{O}}^{\bH_{\mathbb{C}}(\lam)}\otimes_{\mathcal{O}}
F$.

It is clear that the natural homomorphism from
$\ts_F^{\lam}\cong\ts_{\mathcal{O}}^{\lam}\otimes_{\mathcal{O}} F$
to $\td_{\mathcal{O}}^{\lam}\otimes_{\mathcal{O}} F$ is
surjective. Since $\td_F^{\lam}$ is the unique simple head of
 $\ts_F^{\lam}$, it follows that $\td_F^{\lam}$ is also the unique simple head of
$\td_{\mathcal{O}}^{\lam}\otimes_{\mathcal{O}} F$. Hence
$\td_F^{\bH_F(\lam)}\cong\Bigl(\td_F^{\lam}\Bigr)^{\sigma}$ is
also the unique simple head of $
\Bigl(\td_{\mathcal{O}}^{\lam}\Bigr)^{\sigma}\otimes_{\mathcal{O}}
F$. Therefore, $\td_F^{\bH_F(\lam)}$ must also be a composition
factor of
$\td_{\mathcal{O}}^{\bH_{\mathbb{C}}(\lam)}\otimes_{\mathcal{O}}
F$, and hence be a composition factor of
$\ts_F^{\bH_{\mathbb{C}}(\lam)}\cong\ts_{\mathcal{O}}^{\bH_{\mathbb{C}}(\lam)}\otimes_{\mathcal{O}}
F$. Hence $\bH_F(\lam)\trianglelefteq\bH_{\mathbb{C}}(\lam)$.
Using induction on the dominance order $\,\trianglelefteq\,$, it
is easy to see that $\bH_F(\lam)=\bH_{\mathbb{C}}(\lam)$ for any
$\lam\in\mathcal{P}_n$, as required. This completes the proof of
the theorem.  \hfill\qed


\bigskip\bigskip\bigskip
\begin{thebibliography}{99}


\bibitem{A2} S. Ariki, On the decomposition numbers of the
Hecke algebra of $G(m,1,n)$, J. Math. Kyoto Univ. {\bf 36},
789--808 (1996).

\bibitem{A} S. Ariki, On the classification of simple modules
for cyclotomic Hecke algebras of type $G(m,1,n)$ and Kleshchev
multi-partitions, Osaka J. Math. {\bf 38}, (4), 827--837 (2001).

\bibitem{A3} S. Ariki, Proof of the modular branching rule for cyclotomic Hecke algebras,
preprint, math.RT/0511552.

\bibitem{AM} S. Ariki and A. Mathas, The number of simple
modules of the Hecke algebras of type $G(r,1,n)$, Math. Z., {\bf
233}, (3), 601--623 (2000).

\bibitem{AtM} M.F. Atiyah and I.G. Macdonald, Introduction to
Commutative Algebra (Addison-Wesley, 1969).

\bibitem{BO} C. Bessenrodt and J. B. Olsson, Branching of
modular representations of the alternating groups, J. Alg., {\bf
209}, 143--174 (1998).

\bibitem{B} J. Brundan, Modular branching rules and the
Mullineux map for Hecke algebras of type $A$, Proc. London. Math.
Soc., {\bf 77}, (3), 551--581 (1998).

\bibitem{CR} C. W. Curtis and L. Reiner, Methods of
Representations Theory, I (Wiley-Interscience, New York, 1981).

\bibitem{DJ} R. Dipper and G. D. James, Representations of
Hecke algebras of type $B_n$, J. Alg., {\bf 146}, 454--481 (1992).

\bibitem{DJMu} R. Dipper, G. D. James and E. Murphy, Hecke
algebras of type $B_n$ at roots of unity, Proc. London. Math.
Soc., {\bf 70}, (3), 505--528 (1995).

\bibitem{DR} J. Du and H. Rui, Specht modules and branching
rules for Ariki-Koike algebras, Comm. Algebra, {\bf 29},
4701--4719 (2001).

\bibitem{F} B. Ford, Irreducible representations of the alternating group
  in odd characteristic,  Proc. Amer. Math. Soc., {\bf 125}, (2), 375--380 (1997).

\bibitem{FK} B. Ford and A. Kleshchev, A proof of the Mullineux
conjecture, Math. Z., {\bf 226}, 267--308 (1997).

\bibitem{Ge} M. Geck, Representations of Hecke algebras at
roots of unity, S\'eminaire Bourbaki, 50\`eme ann\`ee, 1997--98,
Ast\'erisque, Exp. 836, No. 252, 33--55 (1998).

\bibitem{GR} M. Geck and R. Rouquier, Centers and simple
modules for Iwahori-Hecke algebras, in: Finite reductive groups,
related structures and representations, Progress in Math. Vol.
141, edited by M. Cabanes (Birkh\"auser, Boston, 1997), 251--272.

\bibitem{GL} J. J. Graham and G. I. Lehrer, Cellular algebras,
Invent. Math., {\bf 123}, 1--34 (1996).

\bibitem{G} I. Grojnowski, Affine $sl_p$ controls the
representation theory of the symmetric group and related Hecke
algebras, ArXiv math.RT/9907129.

\bibitem{GV} I. Grojnowski and M. Vazirani, Strong multiplicity
one theorems for affine Hecke algebras of type $A$, Transform.
Groups, {\bf 6}, (2), 143--155 (2001).

\bibitem{Hu1} J. Hu, A Morita equivalence theorem for Hecke
   algebras of type $D_n$ when $n$ is even, Manuscr. Math., {\bf
   108}, 409--430 (2002).

\bibitem{Hu2} J. Hu, Crystal basis and simple modules for Hecke algebra of
  type $D_n$, J. Alg., {\bf 267}, (1), 7--20 (2003).

\bibitem{Hu3} J. Hu, Modular representations of Hecke algebras of type
  $G(p,p,n)$, J. Alg., {\bf 274}, (2), 446--490 (2004).

\bibitem{Hu4} J. Hu, Crystal basis and simple modules for Hecke algebras of
  type $G(p,p,n)$, math.RT/0506555, preprint.

\bibitem{K} M. Kashiwara, On crystal basis, in: Representations of groups,
CMS Conf. Proc., Vol. 16, edited by B. N. Allison and G. H. Cliff
(Amer. Math. Soc. Soc., Providence, RI, 1995), 155--197.

\bibitem{Kl} A. Kleshchev, Branching rules for the modular representations of symmetric groups III; some
  corollaries and a problem of Mullineux, J. London Math. Soc. {\bf 54}, 25--38 (1995).

\bibitem{L} S. Lang, Algebraic Number Theory, second
edition (Springer-Verlag, New York, 1994).

\bibitem{LLT} A. Lascoux, B. Leclerc and J.-Y. Thibon, Hecke
algebras at roots of unity and crystal bases quantum affine
algebras, Comm. Math. Phys., {\bf 181}, 205--263 (1996).

\bibitem{Ma} A. Mathas, The representation theory of Ariki-Koike and cyclotomic $q$-Schur algebras,
  Adv. Studies Pure Math. {\bf 40}, 261--320, Math. Soc. Japan, Tokyo, 2004.

\bibitem{MM} T. Misra, K.C. Miwa, Crystal bases for the basic representations of $U_q(\ksl_n)$,
 Comm. Math. Phys., {\bf 134}, 79--88 (1990).

\bibitem{P} C. Pallikaros, Representations of Hecke algebras of type $D_n$,
  J. Alg., {\bf 169}, 20--48 (1994).

\bibitem{S} J.P. Serre, Local Fields, Graduate texts in
Mathematics Vol.67 (Springer-Verlag, New York, 1979).

\end{thebibliography}
\end{document}